\documentclass[11pt]{amsart}

\usepackage{amsfonts,amssymb,stmaryrd,amscd,amsmath,latexsym,amsbsy}

\usepackage{amssymb}
\usepackage{amsfonts}
\usepackage{latexsym}

\newtheorem{theorem}{Theorem}[section]
\newtheorem{lemma}[theorem]{Lemma}
\newtheorem{proposition}[theorem]{Proposition}
\newtheorem{corollary}[theorem]{Corollary}
\theoremstyle{definition}
\newtheorem{definition}[theorem]{Definition}

\newcommand{\ben}{\begin{enumerate}}
\newcommand{\een}{\end{enumerate}}

\theoremstyle{plain}

\newtheorem*{sol}{Solution}

\theoremstyle{definition}

\theoremstyle{remark}

\newcommand{\solu}[1]{\begin{sol}{\bf (\ref{#1})}}

\begin{document}

\title{New deformations of group algebras of Coxeter groups}

\author{Pavel Etingof and Eric Rains}

\maketitle

\centerline{Dedicated to the memory of Walter Feit}
\vskip .1in

\section{Introduction}

The goal of this paper is to define new deformations of group
algebras of Coxeter groups. Recall that a 
Coxeter group $W$ is generated by elements $s_i, i\in I$
modulo two kinds of relations -- the involutivity relations $s_i^2=1$
and the relations $(s_is_j)^{m_{ij}}=1$, where $2\le m_{ij}=m_{ji}\le
\infty$; in presence of the involutivity relations 
these are equivalent to the {\it braid relations} 
$s_is_js_i...=s_js_is_j...$ ($m_{ij}$ factors). 
The traditional way to deform the group algebra $\Bbb
C[W]$ is to deform the involutivity relations to the Hecke
relations \linebreak $(s_i-q)(s_i+q^{-1})=0$, and keep the braid relations
unchanged. This yields the Hecke algebra 
$\Bbb C_q[W]$ of $W$, which is classically known to be a flat 
deformation of $\Bbb C[W]$. 

On the contrary, the deformation $A$ of $\Bbb C[W]$ 
that we study in this paper is obtained by 
keeping the involutivity relations fixed, and deforming the 
braid relations to 
$$
(s_is_j-t_{ij1})...(s_is_j-t_{ijm_{ij}})=0
$$
when $m_{ij}<\infty$. Here $t_{ijk}$, $k\in \Bbb Z_{m_{ij}}$,
 are new commuting variables 
such that $t_{ijk}=t_{ji,-k}^{-1}$, and
$s_pt_{ijk}=t_{jik}s_p$. We also consider the subalgebra $A_+$
of $A$ generated by $s_is_j$ and $t_{ijk}$. It is a
deformation of the group algebra $\Bbb C[W_+]$ of the group
$W_+$ of even elements of $W$. 
\footnote{In view of the relation $s_pt_{ijk}=t_{jik}s_p$, 
$A$ is not quite an honest deformation of $\Bbb C[W]$ 
(as the variables $t_{ijk}$ are not central). However,
the subalgebra $A_+$ is an honest deformation of 
$\Bbb C[W_+]$, and $A$ is a semidirect product of $\Bbb Z_2$
with $A_+$.}  

A priori, it is not clear that the algebras $A,A_+$ are ``well
behaved''. We show that their ``good'' or ``bad'' behavior is 
determined completely by whether they are formally flat  (i.e. 
whether the corresponding completed algebras over $\Bbb C[[\log
t_{ijk}]]$ are flat deformations of $\Bbb C[W],\Bbb C[W_+]$).
More specifically, we show that if $A$ is formally flat, 
then it is algebraically flat (i.e., free as a module over 
$R:=\Bbb C[t_{ijk}]$), and moreover
any set of reduced words in $s_i$ bijectively representing 
all elements of $W$ defines a basis in $A$ over $R$. 
In particular, this yields a canonical filtration $F^\bullet$ on $A$
(by length of reduced words) such that $F^nA/F^{n-1}A$ is a 
finitely generated free $R$-module. 

Unlike the Hecke algebra $\Bbb C_q[W]$, the algebras $A$,
$A_+$ are not necessarily flat, and we determine when exactly 
this happens. First of all, it is easy to see that flatness always
holds for Coxeter groups of rank 1 and 2. Further, 
we show that for Coxeter groups of rank 3, the 
flatness of $A,A_+$ is equivalent to the condition that 
$W$ is an infinite group. In other words, the (unordered) triple
$m_{12},m_{13},m_{23}$ must be different from 
the triples $(2,2,m)$ ($m<\infty$), $(2,3,3)$, $(2,3,4)$, 
$(2,3,5)$, which correspond to finite Coxeter groups of rank 3: 
$A_1\times I_{2m}$, $A_3,B_3,H_3$. 

This implies that in any
rank, a necessary condition for $A,A_+$ to be flat is that $W$ 
does not contain finite parabolic subgroups of rank 3. 
Our main result asserts that this condition is also sufficient.
The proof is based on consideration of constructible sheaves 
on the cell complex associated to the group $W$. 

The motivation for this paper comes from the fact that for 
rank 3 Coxeter groups $W$, the algebra $A_+$ 
coincides with the Hecke algebra corresponding to the orbifold 
$H/F$, where $H$ is the sphere, Euclidean plane, 
or Lobachevsky plane, and $F$ is the group generated by 
rotations by the angles $2\pi/m_{ij}$
around the vertices of a triangle in $H$ with angles 
$\pi/m_{ij}$. Such Hecke algebras were introduced in \cite{E}, 
and it was shown in \cite{E} (using the theory of Cherednik
algebras and KZ functor) that they are formally flat if $H$ is a
Euclidean or Lobachevsky plane, but not flat if $H$ is a sphere,
which is our main result for rank $3$. (In fact, the proof of
the main result in arbitrary rank is based on similar ideas, 
the difference being that the category of $D$-modules used in
\cite{E} is replaced with the category of constructible sheaves.) 
Thus, as a by-product, we obtain the algebraic PBW theorem for the
Hecke algebras of polygonal Fuchsian groups defined in \cite{E}.

The cases when $H$ is the Euclidean plane and $m_{ij}<\infty$ (i.e., 
$(m_{12},m_{13},m_{23})=(3,3,3),(2,4,4),(2,3,6)$, and
 $W$ is the affine Weyl group of types $A_2$, $B_2$, $G_2$)
were also discussed in \cite{EOR}, as in these cases 
the algebras $A_+$ are the generalized double affine Hecke
algebras of rank 1 of types $E_6,E_7,E_8$ which provide
quantizations of del Pezzo surfaces. Thus we obtain 
new PBW filtrations and bases on these algebras, 
which (unlike those in \cite{EOR}) are constructed without using
a computer. (These filtrations, however, have the flaw that the
idempotent used to define the spherical subalgebra is not
homogeneous). 

{\bf Acknowledgments.} P.E. is very grateful to A. Polishchuk
for an explanation about constructible sheaves, and thanks 
A. Braverman, D. Kazhdan, and C. McMullen for useful
discussions. The work of P.E.
was partially supported by the NSF grant DMS-9988796
and the CRDF grant RM1-2545-MO-03. 
E.R. was supported in part by NSF Grant No.
DMS-0401387.

\section{Definition of the algebras $A(M),A_+(M)$}

\subsection{Coxeter groups}
Recall the basics of the theory of Coxeter groups (see
e.g. \cite{B}).

Let $I$ be a finite set. Let $\Bbb Z_{\ge 2}$ denote the set of
integers which are $\ge 2$. A Coxeter matrix 
over $I$ is a collection $M$ of elements 
$m_{ij}\in \Bbb Z_{\ge 2}\cup \lbrace{\infty\rbrace}$, 
$i,j\in I$, $i\ne j$, such that $m_{ij}=m_{ji}$. 
The rank $r$ of $M$ is, by definition, the cardinality of the set
$I$. 

Let $M$ be a Coxeter matrix. Then one defines the Coxeter group
\footnote{When we talk about Coxeter groups, we always assume
that they are equipped with a fixed system of generators,
which corresponds to the notion of a Coxeter system from \cite{B}.} 
$W(M)$ by generators $s_i$, $i\in I$, and defining relations 
$$
s_i^2=1,\ (s_is_j)^{m_{ij}}=1 \text{ if }m_{ij}\ne \infty.
$$
The group $W(M)$ has a sign character 
$\xi: W(M)\to \lbrace{\pm 1\rbrace}$ given
by $\xi(s_i)=-1$. Denote by $W_+(M)$ the kernel of $\xi$, 
i.e. the subgroup of even elements. 

Let $I'\subset I$ be a subset, and $M'$ the submatrix of $M$ 
consisting of $m_{ij}, i,j\in I'$. Then we have a natural map 
$W(M')\to W(M)$. It is known (\cite{B}) that this map is
injective. Thus $W(M')$ is a subgroup of $W(M)$, which is called
the parabolic subgroup corresponding to $I'$. 

\subsection{Deformed Coxeter group algebras}

Define the algebra $A(M)$ by invertible generators 
$s_i, i\in I$, and $t_{ijk}$, $i,j\in I$, $k\in \Bbb Z_{m_{ij}}$ 
for $(i,j)$ such that $m_{ij}<\infty$ 
with defining relations 
$$
t_{ijk}=t_{ji,-k}^{-1},\
s_i^2=1,
$$
$$
\prod_{k=1}^{m_{ij}}(s_is_j-t_{ijk})=0\text{ if }m_{ij}<\infty, 
$$
$$
[t_{ijk},t_{i'j'k'}]=0,\ 
s_pt_{ijk}=t_{jik}s_p.
$$

Define also the algebra $A_+(M)$ over $R:=\Bbb C[t_{ijk}]$
($t_{ijk}=t_{ji,-k}^{-1}$) by generators 
$a_{ij}$, $i\ne j$ ($a_{ij}=a_{ji}^{-1}$), and relations 
$$
\prod_{k=1}^{m_{ij}}(a_{ij}-t_{ijk})=0\text{ if }m_{ij}<\infty,
$$
$$
a_{ij}a_{jp}a_{pi}=1.
$$

Let $0\in I$ be an element. Then we can define 
an involutive automorphism $\sigma_0$ of $A_+(M)$
(as an algebra over $\Bbb C$) by the formulas  
$\sigma_0(t_{ijk})=t_{jik}$, 
$\sigma_0(a_{ij})=a_{ji}$ if $i=0$ or $j=0$, 
and $\sigma_0(a_{ij})=a_{0i}a_{j0}$ if $i,j\ne 0$. It is easy to show that 
this automorphism is well defined.\footnote{Note that up to inner
automorphisms, $\sigma_0$ does not depend on the choice of 
the element $0\in I$.} 
Thus we can define the semidirect product 
$A_0(M)=\Bbb C\Bbb Z_2\ltimes A_+(M)$ using the automorphism
$\sigma_0$. 

\begin{proposition}\label{iso}
The assignment $f(a_{ij})=s_is_j$, $f(\sigma_0)=s_0$ 
uniquely extends to an isomorphism $f: A_0(M)\to A(M)$. 
\end{proposition}

\begin{proof}
It is easy to check that $f$ uniquely extends to a surjective
homomorphism of algebras. Moreover, 
it is easy to construct the inverse of $f$:
it is given by the formula $f^{-1}(s_i)=\sigma_0a_{0i}$
for $i\ne 0$, and $f(s_0)=\sigma_0$. We are done. 
\end{proof}

Thus, $A(M)=\Bbb C\Bbb Z_2\ltimes A_+(M)$. 

The following proposition explains the connection between 
the algebras $A(M)$, $A_+(M)$ and the groups $W(M)$, $W_+(M)$. 

Let $J$ be the ideal in $R$ generated by the
elements $t_{ijk}-\exp(2\pi {\rm i}k/m_{ij})$. It is easy to see
that $JA(M)=A(M)J$, so $JA(M)$ is a two-sided ideal in $A(M)$. 

\begin{proposition} One has $A(M)/JA(M)=\Bbb C[W(M)]$, 
$A_+(M)/JA_+(M)=\Bbb C[W_+(M)]$. 
\end{proposition}

\begin{proof}
Straightforward.
\end{proof}

\subsection{Spanning sets for $A(M),A_+(M)$}

If $w$ is a word in letters $s_i$, let $T_w$ be the corresponding
element of $A(M)$. 
Choose a function $w(x)$ which attaches
to every element $x\in W(M)$, a reduced word $w(x)$
representing $x$ in $W(M)$. 

We will now prove the following important result.

\begin{theorem}\label{spa} (i) The elements $T_{w(x)}$, $x\in W$, form a spanning set
in $A(M)$ as a left $R$-module.

(ii) The elements $T_{w(x)}$, $x\in W_+$, form a spanning set
in $A_+(M)$ as a left $R$-module.
\end{theorem}

\begin{proof} It is clear that (ii) follows from (i), so it
suffices to prove (i). 

Let us write the relation
$$
\prod_{k=1}^{m_{ij}}(s_is_j-t_{ijk})=0
$$
as a deformed braid relation:
$$
s_js_is_j...+S.L.T.=t_{ij}s_is_js_i...+S.L.T., 
$$
where $t_{ij}=(-1)^{m_{ij}+1}t_{ij1}...t_{ijm_{ij}}$, S.L.T.
mean ``smaller length terms'', 
and the products on both sides have length $m_{ij}$.
This can be done by multiplying the relation by $s_is_j...$ 
($m_{ij}$ factors).

Now let us show that $T_{w(x)}$ span $A(M)$ over $R$.
Clearly, $T_w$ for all words $w$ span $A(M)$. 
So we just need to take any word $w$ and
express $T_w$ via $T_{w(x)}$. 

It is well known from the theory of Coxeter groups (see e.g. \cite{B}) that 
using the braid relations, one can turn any non-reduced word
into a word that is not square free, and any
reduced expression of a given element 
of $W(M)$ into any other reduced expression of the same element.   
Thus, if $w$ is non-reduced, then by using the deformed braid relations
we can reduce $T_w$ to a linear combination of $T_u$ with words $u$ of
smaller length than $w$. 
On the other hand, if $w$ is a reduced expression for some
element $x\in W$, then using the deformed braid relations 
we can reduce $T_w$ to
a linear combination of $T_u$ with $u$ shorter than $w$, and
$T_{w(x)}$.
Thus $T_{w(x)}$ are a spanning set.
The theorem is proved. 
\end{proof}

Thus, $A_+(M)$ is a ``deformation'' of $\Bbb C[W_+(M)]$ over $R$, and  
similarly $A(M)$ is a ``twisted deformation'' of $\Bbb C[W(M)]$.

\section{Flat Coxeter matrices}

\subsection{Definition of a flat Coxeter matrix}

Denote by $\hat A(M)$, $\hat A_+(M)$ the formal versions 
of $A(M)$, $A_+(M)$, i.e., algebras generated (topologically) 
by the same generators
and relations, but with $t_{ijk}=e^{2\pi {\rm i}k/m_{ij}}e^{\tau_{ijk}}$, where
$\tau_{ijk}$ are formal parameters. By virtue of 
Theorem \ref{spa}, $\hat A_+(M)$ is an algebra over $\hat R:=\Bbb
C[[\tau_{ijk}]]$ which is a formal deformation 
of $\Bbb C[W_+]$ with deformation parameters
$\tau_{ijk}$. 

\begin{definition} 
We say that $M$ is a flat Coxeter matrix if 
$\hat A_+(M)$ is a flat deformation of $\Bbb C[W_+]$,
i.e. if $\hat A_+(M)$ is a topologically free left $\hat R$-module.  
\end{definition}

Since $\hat A(M)=\Bbb C\Bbb Z_2\ltimes \hat A_+(M)$, 
for a flat Coxeter matrix we also have 
that $\hat A(M)$ is topologically free as a left $\hat R$-module.

\subsection{Bases of $A(M),A_+(M)$ for flat $M$}

\begin{proposition} Let $M$ be a flat Coxeter matrix. Then

(i) The elements $T_{w(x)}$, $x\in W$, form a basis
in $A(M)$ as a left $R$-module.

(ii) The elements $T_{w(x)}$, $x\in W_+$, form a basis
in $A_+(M)$ as a left $R$-module.
\end{proposition}

\begin{proof}
It is sufficient to show that $T_{w(x)}$ are linearly
independent. This follows from the fact that they are linearly
independent in $\hat A(M)$, which is a consequence of the flatness of
$M$. 
\end{proof}

\begin{corollary}
The $R$-modules $A(M)$ and $A_+(M)$ are free. 
Moreover, they carry a filtration $F^\bullet$,
defined by the condition that $F^n A(M)$ is spanned by
$T_{w(x)}$ for $x$ of length $\le n$. This filtration has the
property that the $R$-modules ${\rm gr}_n A(M)$ and ${\rm gr}_n
A_+(M)$ are finitely generated and free.   
\end{corollary}

We note that the filtration $F^\bullet$ is canonical, i.e., 
independent on the choice of the function $w(x)$. 

{\bf Remark.} Let $\Gamma(M)$ be the graph whose vertices are
elements of $I$, and $i,j$ are connected if $m_{ij}$ is odd. 
Then the filtration $F^n$ can be refined to 
a multi-filtration $F^{\bold n}$, where 
$\bold n$ is a nonnegative integer function on $I$ 
which is constant on the connected components of $\Gamma(M)$. 
Namely, $F^{\bold n}A(M)$ is spanned by $T_{w(x)}$ with $w(x)$
involving $\le \bold n(i)$ copies of $s_i$ for each $i$. 
It is easy to see from the above that 
${\rm gr}_{\bold n} A(M)$ is finitely generated and free
over $R$. 

\subsection{When is a Coxeter matrix flat?}

Let us study the question when a given Coxeter matrix is flat
(i.e., when the algebras $A(M),A_+(M)$ are ``meaningful'').

The case of rank 1 is trivial. In rank 2, the algebra 
$A_+(M)$ is $\Bbb C[a,a^{-1}]$ if $m_{12}=\infty$ and 
$\Bbb C[a,a^{-1}]/(P(a))$ where
$P(a)=(a-t_{121})...(a-t_{12m_{12}})$ otherwise. Thus, any Coxeter matrix of
rank 2 is flat. 

However, in rank 3, we have a much more interesting situation. 
Namely, we have the following theorem. 

\begin{theorem}\label{r3} A Coxeter matrix $M$ of rank 3 is flat 
if and only if the group $W(M)$ is infinite. 
\end{theorem}

Coxeter matrices of rank 3 are conveniently written as 
triples of numbers $m_{12},m_{23},m_{31}$ (the order
does not matter). 
Recall that the Coxeter matrices of rank 3 producing a finite Coxeter group are 
the following: 

1. $M=(2,2,m)$, $m<\infty$ (type $A_1\times I_{2m}$). 

2. $M=(2,3,3)$ (type $A_3$)

3. $M=(2,3,4)$ (type $B_3$)

4. $M=(2,3,5)$ (type $H_3$)

Thus, the theorem claims that a Coxeter matrix is flat if and only
if it does not fall into the four cases listed above.

\begin{proof}
{\bf If.} Let $M=(m_{12},m_{23},m_{31})$. The algebra $A_+(M)$ 
is defined over $R$ by  generators $a_{12},a_{23},a_{31}$ 
with defining relations 
$$
(a_{ij}-t_{ij1})...(a_{ij}-t_{ijm_{ij}})=0,\ ij=12,23,31,
$$
and 
$$
a_{12}a_{23}a_{31}=1.
$$
This algebra is a deformation of the group algebra of the group 
$F=F_{m_{12},m_{23},m_{31}}$ generated by 
$a_{12},a_{23},a_{31}$ with defining relations
$a_{ij}^{m_{ij}}=1$,
$a_{12}a_{23}a_{31}=1$. The group $F$ is isomorphic to the group   
generated by rotations by angles $2\pi/m_{ij}$ 
around vertices of a triangle whose angles 
are $\pi/m_{ij}$. This triangle lies on the sphere, plane, and
Lobachevsky plane if the quantity
$S=\frac{1}{m_{12}}+\frac{1}{m_{23}}+
\frac{1}{m_{31}}$ is $>1,=1$, and $<1$, respectively. 
The cases 1,2,3,4 when $W(M)$ is finite correspond 
exactly to the case $S>1$. Now, the flatness of the algebra 
$\hat A_+(M)$ for $S\le 1$ follows from Theorem 3.3 of \cite{E}
(as $\hat A_+(M)$ is the Hecke algebra of the orbifold 
$H/F$, where $H$ is the plane or Lobachevsky plane,
depending on whether $S=1$ or $S<1$). Another proof of flatness 
of $\hat A_+(M)$ for $S=1$ is given in \cite{EOR} 
(in this case $\hat A_+(M)$ is a generalized double affine Hecke algebra). 
This proves the ``if'' direction of the theorem. 

{\bf Only if.} Suppose $W(M)$ is finite. Assume the contrary,
i.e., that $M$ is flat. Then the algebra $\hat A_+(M)$ 
is a free module over $\hat R$ of dimension $D=|W_+(M)|$.
The eigenvalues of $a_{ij}$ in the regular representation  
are equal to $t_{ijk}$, and occur with multiplicity $D/m_{ij}$. 
Thus in the regular representation we have 
$\det(a_{ij})=(\prod_k t_{ijk})^{D/m_{ij}}$. 
So taking the determinant of the relation $a_{12}a_{23}a_{31}=1$, 
we have 
$$
(\prod_{k=1}^{m_{12}} t_{12k})^{D/m_{12}}
(\prod_{k=1}^{m_{23}} t_{23k})^{D/m_{23}}
(\prod_{k=1}^{m_{31}} t_{31k})^{D/m_{31}}=1.
$$
This is a nontrivial relation on $t_{ijk}$, which contradicts to
the flatness of $M$. 
The theorem is proved. 
\end{proof}

Now we are ready to state our main result. 

\begin{theorem}\label{mainc} A Coxeter matrix $M$ is flat if and only if 
for any 3-element subset $\lbrace{i,j,k\rbrace}$ of $I$, 
the Coxeter group generated by $s_i,s_j,s_k$ is infinite. 
\end{theorem}

The theorem is proved in the next susbection. 
Note that in the process of proving Theorem \ref{mainc}, we give
a new proof of the ``if'' part of Theorem \ref{r3}, 
which relies on arguments from topology 
(constructible sheaves) rather than complex
analysis (D-modules used in \cite{E}).

\subsection{Proof of Theorem \ref{mainc}}

Let us first prove the easier ``only if'' direction.
Let $I'=\lbrace{i,j,k\rbrace}\subset I$. 
Let $M'$ be the corresponding Coxeter matrix.
Since $W(M')\subset W(M)$, 
the flatness of $A_+(M)$ implies the flatness of $A_+(M')$. 
But if $A_+(M')$ is flat then by Theorem \ref{r3}, 
$W(M')$ is an infinite group, as desired. 

Now let us prove the ``if'' direction. To do so, 
we introduce a 2-dimensional cell complex $\Sigma$ 
attached to $M$ as follows. The zero-dimensional cells 
of $\Sigma$ are the elements of $W:=W(M)$. 
The 1-dimensional cells are edges connecting 
$w$ and $s_iw$ for each $i\in I$. 
Then we have cycles $w,s_iw,s_js_iw,s_is_js_iw...,s_jw$ of length
$2m_{ij}$ (if $m_{ij}<\infty$), and the 2-dimensional cells of $\Sigma$ 
are $2m_{ij}$-gons attached to these cycles. 

It is easy to see that $W$ acts properly discontinuously on
$\Sigma$, hence so does its subgroup $W_+$. Moreover, it is clear
that the only fixed points of the $W_+$-action on $\Sigma$ are
the  centers of the $2m_{ij}$-gons, with stabilizer 
$\Bbb Z_{m_{ij}}$. Thus we can define an orbifold cell complex 
$Y:=\Sigma/W_+$. It has two vertices (the north and south pole, $N,S$), edges 
$e_i$ between them corresponding to $s_i$, and for each 2-element
set $\lbrace{i,j\rbrace}\subset I$, a 
disk $D_{ij}$ whose boundary is identified with the circle made up by
$e_i$ and $e_j$. The disk $D_{ij}$ has an orbifold point $z_{ij}$
in the center, whose isotropy group is $\Bbb Z_{m_{ij}}$. 

We will need the following theorem (see e.g. \cite{DM}). 

\begin{theorem}\label{dav}
If $W$ has no finite 
parabolic subgroups of rank 3 then $\Sigma$ is contractible.
\end{theorem}

Let ${\mathcal C}$ be the category of constructible sheaves 
of complex vector spaces on 
$\Sigma$ with respect to the stratification into cells
(see \cite{Sch} for definitions). 
Let $\Bbb C$ be the constant sheaf on $\Sigma$.

\begin{lemma}\label{po} For any Coxeter matrix we have 
${\rm Ext}^1_{\mathcal C}(\Bbb C,\Bbb C)=H^1(\Sigma, \Bbb
C)$, and ${\rm Ext}^2_{\mathcal C}(\Bbb C,\Bbb C)\subset
H^2(\Sigma,\Bbb C)$. 
\end{lemma}

\begin{proof}\footnote{This argument was provided to us by A. Polishchuk.}
Let $\widetilde{\mathcal C}$ be the abelian category
of all sheaves of complex vector spaces on $\Sigma$. Then ${\mathcal C}$ is 
full abelian subcategory closed under extensions.
In this situation for any 
$F,G\in {\mathcal C}$,  
${\rm Ext}^1_{\mathcal C}(F,G)={\rm Ext}^1_{\widetilde {\mathcal C}}(F,G)$,
and the natural map 
${\rm Ext}^2_{{\mathcal C}}(F,G)\to {\rm Ext}^2_{\widetilde{\mathcal
C}}(F,G)$ is injective. But it is well known that 
${\rm Ext}^i_{\widetilde{\mathcal C}}(\Bbb C,\Bbb C)=H^i(\Sigma,\Bbb C)$. 
This implies the statement. 
\end{proof}

{\bf Remark.} As was explained to us by A. Polishchuk, 
the inclusion 
${\rm Ext}^2_{\mathcal C}(F,G)\to
{\rm Ext}^2_{\widetilde{\mathcal C}}(F,G)$ is in fact an isomorphism. 
It follows from the fact
that our stratification satisfies the property that the closure
of every cell is homeomorphic to a closed ball (in a way
compatible with boundaries). Indeed, it suffices to prove this
when $F$ and $G$ are simple.
Then it follows from the fact that the algebra of
${\rm Ext}^*_{\widetilde{\mathcal C}}$
between simple objects of ${\mathcal C}$ 
is generated in degree $1$
(apply Corollary 2.2 of \cite{P} for zero perversity).

\begin{lemma}\label{ext=0} If $W$ has no finite parabolic
subgroups of rank $3$, one has 
${\rm Ext}^j_{\mathcal C}(\Bbb C,\Bbb C)=0$ for $j=1,2$.
\end{lemma}

\begin{proof}
The result follows from Theorem \ref{dav} and Lemma \ref{po}.
\end{proof} 

Let ${\mathcal D}$ be the category of constructible sheaves 
on the orbifold $Y$ with respect to stratification into cells. 
Thus, an object of ${\mathcal D}$ is a 
constructible sheaf on the complement $Y'\subset Y$ 
of the points $z_{ij}$ with respect to the same
stratification (intersected with $Y'$), such that the monodromy $g_{ij}$ 
around $z_{ij}$ satisfies the equation $g_{ij}^{m_{ij}}=1$. 
 
Let $\pi: \Sigma\to
Y$ be the natural projection, 
$\pi_!: {\mathcal C}\to {\mathcal D}$ be the direct image
functor with compact supports, and ${\bold M}=\pi_!(\Bbb C)$. Thus ${\bold M}$ is a 
local system on the orbifold $Y$. 
The monodromy representation of ${\bold M}$ over $Y'$
is the regular representation of $W_+$. 

\begin{lemma}\label{cohom}
One has ${\rm Ext}^j_{\mathcal D}({\bold M},{\bold M})=0$ for $j=1,2$. 
\end{lemma}

\begin{proof}
$$
{\rm Ext}^j_{{\mathcal D}}({\bold M},{\bold M})=
{\rm Ext}^j_{{\mathcal D}}(\pi_!\Bbb C,\pi_!\Bbb C)=
$$
$$
{\rm Ext}^j_{{\mathcal C}}(\Bbb C,\pi^*\pi_!\Bbb C)=\Bbb C[W_+]\otimes 
{\rm Ext}^j_{{\mathcal C}}(\Bbb C,\Bbb C).
$$
But ${\rm Ext}^j_{{\mathcal C}}(\Bbb C,\Bbb C)=0$ by Lemma \ref{ext=0}. 
We are done. 
\end{proof} 

A basic fact about constructible sheaves,
going back to Fulton, Goresky, McCrory, and MacPherson
(\cite{Sh}; see also \cite{Vy}),
is that the category of constructible
sheaves on a cell complex with respect to the stratification
into cells is equivalent to the category of ``cellular sheaves'',
i.e., modules over a certain algebra $B$
(path algebra of a certain quiver with relations). 

Let us construct $B$ in the case of the orbifold cell complex
$Y$. Let us fix an ordering on the set $I$. 
Let ${\mathcal S}\in {\mathcal D}$. Let $V_N$, $V_S$, $V_i$,
$V_{ij}$ ($i\ne j$) be the
spaces of sections of ${\mathcal S}$ over 
small neighborhoods of $N,S$, the
midpoints of $e_i$, and the points of $D_{ij}$ close to the midpoints
of $e_i$, respectively. We then have natural restriction maps
$f_{Ni}: V_N\to V_i$, $f_{Si}: V_S\to V_i$, 
$h_{ij}: V_i\to V_{ij}$, and monodromy maps 
$g_{ij}: V_{ij}\to V_{ji}$, with relations 
$$
g_{ij}h_{ij}f_{Ni}=g_{ji}f_{Nj},\ \text{ for }i<j,
$$
$$
g_{ij}h_{ij}f_{Si}=g_{ji}f_{Sj},\ \text{ for }i>j,
$$
and 
$$
 (g_{ij}g_{ji})^{m_{ij}}=1.
$$
Thus, ${\mathcal S}$ defines a representation of a certain quiver
with relations. We define the algebra $B$ to be the path 
algebra of this quiver (modulo the relations). It is 
then well known that category ${\mathcal D}$ is equivalent to the
category of $B$-modules: the equivalence sends ${\mathcal S}$ to
the $B$-module 
$$
V:=V_N\oplus V_S\oplus\oplus_i V_i\oplus \oplus_{i,j}V_{ij}.
$$  
 
Now let $\tau=(\tau_{ijk}), k\in \Bbb Z_{m_{ij}}$, 
be a collection of formal parameters.
Define the algebra $B(\tau)$ by the same generators as 
$B$, and the same relations except for one modification: 
the relation $(g_{ij}g_{ji})^{m_{ij}}=1$ 
is replaced with 
$$
(g-t_{ij1})...(g-t_{ijm_{ij}})=0,
$$
where $g:=g_{ij}g_{ji}$. 
(we recall that $t_{ijk}:=e^{2\pi {\rm
i}k/m_{ij}}e^{\tau_{ijk}}$). It is clear that
$B(\tau)/(\tau=0)=B$. 

We will need the following flatness result. 

\begin{proposition}\label{fl} 
For any Coxeter matrix $M$, the algebra $B(\tau)$ is a flat
deformation of $B$. 
\end{proposition}

\begin{proof} For any $i<j$ let $p_{ij}$ be the idempotent of $B$ 
which acts by $1$ on $V_N,V_S,V_i,V_j$ and on $V_{ij},V_{ji}$ if 
$m_{ij}<\infty$, and acts by $0$
on all the other spaces. Then the direct sum of right modules 
$p_{ij}B$ over $B$ is faithful, so it suffices to show that they
can be deformed to right $B(\tau)$-modules. 
Clearly, $B$ preserves the kernel of $p_{ij}$, so
$p_{ij}BB=p_{ij}Bp_{ij}$, and thus $B_{ij}:=p_{ij}B$ is a unital algebra
with unit $p_{ij}$. Replacing relations as above, we can define 
a deformation $B_{ij}(\tau)=p_{ij}B(\tau)p_{ij}$ of $B_{ij}$, 
and it suffices to show
that this deformation is flat. 

We clearly only need to consider the case $m_{ij}<\infty$. 
Let us consider the category ${\mathcal E}_{ij}$ 
of representations of $B_{ij}$ in which all the 
quiver arrows are isomorphisms. It is easy to see 
(by explicitly writing a basis of $B_{ij}$) that the direct sum
of these representations is faithful. Thus it suffices 
to show that any such representation can be deformed to a
representation of $B_{ij}(\tau)$. But it is easy to see 
that the category ${\mathcal E}_{ij}$ is equivalent to the category 
of vector spaces $U$ with a linear map $g: U\to U$ 
such that $g^{m_{ij}}=1$. This shows that any object
of this category can be easily deformed to a module over
$B_{ij}(\tau)$ (by deforming $g$ to an operator 
satisfying the equation $(g-t_{ij1})...(g-t_{ijm_{ij}})=0$), as desired.
\end{proof}

Now we can finish the proof of Theorem \ref{mainc}. 
We can regard ${\bold M}$ as 
a $B$-module (in which all the arrows are isomorphisms). 
By Lemma \ref{cohom}, ${\rm Ext}^1_B({\bold M},{\bold M})={\rm Ext}^2_B({\bold M},{\bold M})=0$. 
This implies that ${\bold M}$ can be uniquely deformed to a module
${\bold M}_\tau$ over $B_\tau$. The module ${\bold M}_\tau$ is a quiver
representation where all arrows are isomorphisms, 
so it represents a local system on $Y'$ (over $\Bbb C[[\tau]]$).
The monodromy of this local system is a representation 
of $A_+(M)$, deforming the regular representation of $\Bbb
C[W_+]$. The existence of such deformation implies the flatness
of $A_+(M)$. Theorem \ref{mainc} is proved. 

\subsection{Flatness of Hecke algebras of 
polygonal Fuchsian groups}

Let $\Gamma=\Gamma(m_1,...,m_r)$, $r\ge 3$, be the Fuchsian group 
defined by generators $c_j$, 
$j=1,...,r$, with defining relations 
$$
c_j^{m_j}=1, 
\ \prod_{j=1}^r c_j=1.
$$
Here $2\le m_j\le \infty$. 

In \cite{E} the first author defined the Hecke algebra 
of $\Gamma$, ${\mathcal H}(\Gamma)$, 
by the same (invertible) generators and 
relations 
$$
(c_j-t_{j1})...(c_j-t_{jm_j})=0,\text{ if }m_j<\infty,\ 
\prod_{j=1}^r c_j=1,
$$
where $t_{jk}$ are invertible variables.

\begin{theorem} 
The algebra ${\mathcal H}(\Gamma)$ is free as a left module 
over $R:=\Bbb C[t,t^{-1}]$
if and only if $\sum_j (1-1/m_j)\ge 2$ (i.e., $\Gamma$ is
Euclidean or hyperbolic). 
\end{theorem}

The formal version of this theorem is proved in \cite{E}. 

\begin{proof} Let $M$ be the Coxeter matrix of rank $r$ 
such that $m_{i,i+1}:=m_i$ 
for $i\in \Bbb Z_r$, and $m_{ij}=\infty$ otherwise. 
It is easy to deduce from Theorem \ref{mainc} that
$M$ is flat if and only if 
$\sum_j (1-1/m_j)\ge 2$.   
But for such matrix $A_+(M)={\mathcal
H}(\Gamma)$. We are done. 
\end{proof}

Note that since ${\mathcal H}(\Gamma)=A_+(M)$, we actually obtain
a basis of ${\mathcal H}(\Gamma)$, given by $T_{w(x)}$ for even
$x$, and also a canonical filtration on ${\mathcal H}(\Gamma)$
with free finitely generated quotients.

\end{document}